\documentclass[12pt]{amsart}
\usepackage{amsfonts,amssymb,latexsym,eufrak}

\def\det{\mathrm{det}} 
 
\def\corank{\mathrm{corank}}

\def\Gr{\mathrm{Gr}} 
 
\def\GL{\mathrm{GL}} 
\def\iso{\cong}
\def\isom{\cong}
 
\def\restr{\!\upharpoonright\!}
 
\def\phi{\varphi} 
\def\epsilon{\varepsilon} 
\def\theta{\vartheta} 
 
\def\Z{\mathbb{Z}} 
\def\N{\mathbb{N}} 
\def\C{\mathbb{C}} 

\def\d{\mathrm{d}}

\def\c{\complement}
\def\phtilde{\widetilde{\phantom{y}}}
\def\pt{\mathrm{pt}}

\def\ts{\mathrm{Ts}} 
\def\Diff{\mathrm{Diff}} 
\def\Hom{\mathrm{Hom}}
 
\def\Sym{\mathrm{Sym}}

\def\codim{\mathrm{codim}} 
\def\im{\mathrm{Im}} 
\def\ker{\mathrm{Ker}} 
\def\coker{\mathrm{Coker}}

\newcommand{\blds}[1]{\mbox{\boldmath $\scriptstyle #1 $}}
\newcommand{\gbinom}[2]{ {{#1}\brace{#2}} }

\theoremstyle{plain} 

\newtheorem{fact}{Fact}[section] 
\newtheorem*{conj}{Conjecture} 
\newtheorem{corr}[fact]{Corollary} 
\newtheorem{lem}[fact]{Lemma} 
\newtheorem{prop}[fact]{Proposition} 
\newtheorem{thm}[fact]{Theorem}
 
\theoremstyle{definition} 

\newtheorem*{rem}{Remark} 
\newtheorem{defn}[fact]{Definition} 
\newtheorem*{exmp}{Example}

\theoremstyle{plain} 

\begin{document}
\title{On second order Thom-Boardman singularities}
 
\subjclass[2000]{Primary 32S20; secondary 57R45.} 
 
\keywords{Thom polynomials, singularities, classes of degeneracy loci, 
global analysis.} 
 
\author[L. M. Feh\'er]{L\'aszl\'o M. Feh\'er} 
\address{Department of Analysis, E\"otv\"os University, Budapest} 
\email{lfeher@cs.elte.hu} 
 
\author[B. K\H om\H uves]{Bal\'azs K\H{o}m\H{u}ves}
\address{Department of Mathematics, Central European University, Budapest} 
\email{komuves@renyi.hu} 

\begin{abstract} 
In this paper we derive closed formulas for the Thom polynomials of two
families of second order Thom-Boardman singularities $\Sigma^{i,j}$.
The formulas are given as linear combinations of Schur polynomials, and all 
coefficients are nonnegative. 
\end{abstract} 
        
\thanks{Both authors are partially supported by the OTKA grant T046365MAT} 
 
\maketitle 

%--------------------------------------------------------------------
\section{Introduction} 
 In the fifties, R.~Thom \cite{thom} defined for generic smooth 
(resp.~analytic) maps $f:N^n\to P^p$ between real (resp.~complex) manifolds, 
singularity subsets $\Sigma^I(f)\subset N$ for each  partition 
$I=(i_1,\dots,i_k)$ inductively, by setting  
\begin{eqnarray*} 
\Sigma^i(f) &= &\{x\in N\::\: \corank(d_xf)=i\} \\ 
\Sigma^{I,j}(f) &=&\Sigma^{j}(f\restr \Sigma^I(f)) 
\end{eqnarray*} 
Later J.~M.~Boardman \cite{boardman} gave a more precise 
and also more  general definition in terms of jet bundles. 
Thom also observed that for generic maps, the cohomology
classes $[\Sigma^I(f)]\in H^*(N)$  represented by (the closures of) these 
subsets depend only on the Stiefel-Whitney (resp.~Chern) classes of the bundles 
$TN$ and $f^*TP$, and that they are given by a universal polynomial 
$[\Sigma^I(n,p)](c(TN),c(f^*TP))$ in these classes, now called the Thom 
polynomial of the singularity.  
In fact, the same holds for many other singularity classes, that is, for 
(the closure of) a submanifold of some jet space invariant to the 
left-right action of jets of (germs of) biholomorphisms; and for 
\emph{stable} classes of singularities, these polynomials depend 
only the Chern classes of the virtual normal bundle $f^*TP-TN\in K(N)$ 
and the codimension $r=p-n$ (see Proposition \ref{stable}). Since the Thom-Boardman singularities are
stable, their Thom polynomials can be expressed as a polynomial 
$[\Sigma^{I}(r)]$ in the formal variables $c_1,c_2,c_3,\dots$.
 
So far, very few of these polynomials are known explicitly. 
The Thom polynomials for $\Sigma^i$ were calculated by Porteous 
\cite{porteous}; Ronga \cite{ronga} gave an algorithm to calculate 
the classes $\Sigma^{i,j}$, but this algorithm is inefficient even 
using today's top personal computers, and the only known formula 
derived from it gives the Thom polynomials of $\Sigma^{1,1}$ (for all $r$). 
M. Kazarian in \cite{kaza:localized} simplified the argument of Ronga 
and also gave slightly different versions of the original algorithm. Recently the 
Thom polynomials for $\Sigma^{1,1,1}$ was calculated \cite{bfr} using the
methods of restriction equations; and some other sporadic 
results are known (for example $\Sigma^{1,1,1,1}$ in the case $r=0$, 
see \cite{gaffney};  $\Sigma^{1^k}$ for $k\le 8$, $r=0$, see \cite{rrthompol}). 
The present work can be considered as the first step of the more ambitious 
quest for a \emph{formula} for the Thom polynomials of \emph{all} second order 
Thom-Boardman singularities $\Sigma^{i,j}$. We would like to emphasize that it 
wasn't a priori clear that such formulas exist: In fact, the method of 
restriction equations \cite{rrthompol} suggested 
just the contrary, since the 
complexity of singularities having smaller codimension than $\Sigma^{i,j}(r)$ 
is increasing rapidly with $i$.

The success of the method of restriction equations (for a more 
detailed account than in \cite{rrthompol} see \cite{fr}) encouraged 
us to try it in this situation as well. To our surprise we realized 
that the restriction equation we started to study is a consequence 
of the fact mentioned above, that these singularities are stable. 
It turned out that in essence, everything needed to produce our
formulas was known for at least 20 years.

From the results of Porteous \cite{porteous} and Ronga \cite{ronga} one 
can deduce that in some sense it is natural to write the polynomials 
$[\Sigma^{i,j}(r)]$ in terms of Schur polynomials. After finishing this work
we were informed by Piotr Pragacz that he also used Schur polynomial methods
to calculate other Thom polynomials \cite{pragacz-thom}.
Theorems \ref{main1} and \ref{main2nice} 
give the formulas as (nonnegative) 
linear combinations of Schur polynomials. To state them, introduce the
notations 
\[ E_{\lambda/\mu}(n):=\det \left[ \binom{\lambda_k+n-k}{\mu_l+n-l} \right ]_{n\times n} 
\quad\quad\textrm{and}\quad\quad
F_{\lambda/\mu}(n):=\det \left[ \gbinom{\lambda_k+n-k}{\mu_l+n-l} \right ]_{n\times n}
\]
where $\lambda$ and $\mu$ are partitions, and $\gbinom{n}{k}=\sum_{j=0}^k\binom{n}{j}$.

\theoremstyle{plain} 
\newtheorem*{thm_main1}{Theorem \ref{main1}} 
\begin{thm_main1}
The Thom polynomial of the second order Thom-Boardman singularity $\Sigma^{i,j}(-i+1)$ is 
\[ [\Sigma^{i,j}(-i+1)] =  \sum_{\mu\subset\delta} 2^{|\mu|-j(j-1)/2} \cdot E_{\delta/\mu}(i) \cdot s_{({d-|\mu|},\widetilde\mu)} \]
where $\delta$ is the `staircase' partition $\delta=(j,j-1,\dots,2,1)$,
and $d$ is the codimension of $\Sigma^{i,j}(-i+1)$, that is, $d=i+j(j+1)/2$.
\end{thm_main1}

\noindent
Note that these are the simplest singularities which can occur for negative codimension maps.
In the setting $f:N^n\to P^{n-i+1}$ described above, $s_\lambda$ becomes $s_\lambda(f^*TP-TN)$.

\newtheorem*{thm_main2nice}{Theorem \ref{main2nice}} 
\begin{thm_main2nice}
The Thom polynomial of $\Sigma^{i,1}(r)$ is 
\[ {}[\Sigma^{i,1}(r)] = 
\sum_{(\nu,\mu)\in I} F_{\nu/\mu}(i)\cdot
 s_{(i^{(r+i)}+\c\widetilde\nu,\widetilde\mu)} 
\]
where $I=\{(\nu,\mu )\::\:\nu\subset (r+i)^i,\; l(\mu)\le i,\; |\nu|-|\mu|=i-1\}$ (see Figure~\ref{abra}).
\end{thm_main2nice}

\begin{figure}\label{abra}
\setlength{\unitlength}{1mm}
\begin{picture}(135,80)(-10,-5)
\linethickness{0.15mm}
\put(0,0){\line(1,0){120}}
\put(0,0){\line(0,1){65}}
\put(50,0){\line(0,1){60}}
\put(0,30){\line(1,0){60}}
\put(0,60){\line(1,0){50}}
\linethickness{1mm}
\put(-0.5,0){\line(1,0){111}}
\put(0,-0.5){\line(0,1){61}}
\put(-0.5,60){\line(1,0){31}}
\put(30,60.5){\line(0,-1){6}}
\put(29.5,55){\line(1,0){6}}
\put(35,55.5){\line(0,-1){16}}
\put(34.5,40){\line(1,0){6}}
\put(40,40.5){\line(0,-1){6}}
\put(39.5,35){\line(1,0){11}}
\put(50,35.5){\line(0,-1){6}}
\put(49.5,30){\line(1,0){20.9}} %hack
\put(70,30.5){\line(0,-1){11}}
\put(69.5,20){\line(1,0){6}}
\put(75,20.5){\line(0,-1){6}}
\put(74.5,15){\line(1,0){16}}
\put(90,15.5){\line(0,-1){11}}
\put(89.5,5){\line(1,0){21}}
\put(110,5.5){\line(0,-1){6}}
\put(-3,28.5){$i$}
\put(-5,58.5){$2i$}
\put(45,-4){$h=r+i$}
\put(21,14){$i^{(r+i)}$}
\put(17,44){$\c\widetilde\nu$}
\put(60,12){$\widetilde\mu$}
\put(42,48){$\widetilde\nu$}
\end{picture}
\caption{The partition $(i^{(r+i)}+\c\widetilde\nu,\widetilde\mu)$.}
\end{figure}

The second order Thom-Boardman singularities are indexed by the three parameters $i,j$ 
and $r$ (although to us, it seems to be more natural to use $h=r+i$ instead of $r$). 
Both Theorem \ref{main1} and \ref{main2nice} provides a closed formula for a 
two parameter family. These provide two ``transversal'' planes in the $3$-dimensional 
space of parameters. In the latter case we get a family where $r$ can run independently, 
giving infinitely many examples of {\it Thom series}, see \cite{fr-d-stab}.
Namely, if we fix $i$ and $j$, the Thom polynomials for different $r$'s should fit together into one 
series
\[ \ts(\Sigma^{i,j} ) = \sum_{\gamma} c_\gamma \cdot d_\gamma \]
where $\gamma$ runs over nonincreasing $\Z$-valued sequences of length $i+k=i+ij-\binom{j}{2}$
with fixed sum $|\gamma| = j(j-i)\le 0$. Here, $c_\gamma$ are (integer)
coefficients---\emph{independent of} $r$---and $d_\gamma$ are ``renormalized'' Schur polynomials
$d_\gamma = s_{ (h^{(i+k)}+\gamma)\phtilde}$.
Furthermore, Theorem \ref{rongacor} says that $c_\gamma=0$ unless 
$\gamma_l\ge 0$ for $l\le i$ and $\gamma_l\le 0$ for $l>i$.
In this notation, Theorem \ref{main2nice} becomes
\[\ts(\Sigma^{i,1})=\sum_{(\nu,\mu)\in I'} F_{\nu/\mu}(i)\cdot d_{(\mu,-\nu_i,-\nu_{i-1},\dots,-\nu_1)}\]
with
$I'=\{(\nu,\mu )\::\: l(\nu)\le i,\; l(\mu)\le i,\; |\nu|-|\mu|=i-1\}$. 
Theorem \ref{main1} provides the ``lowest terms'' for the Thom series  $\ts(\Sigma^{i,j})$.\\
  
We thank Anders Buch, Maxim Kazarian and Rich\'ard Rim\'anyi for conversations on the topic.
We used John Stembridge's SF package for Maple \cite{SF} extensively during the preparation 
of this paper.\\

\subsection*{Notations}
A partition is a nonincreasing sequence of positive 
integers $\mu=(\mu_1\ge\mu_2\ge\dots\ge\mu_n>0)$; the length of a partition 
is denoted by $l(\mu)=n$, its weight by $|\mu|=\sum \mu_i$. We adopt the 
convention that $\mu_i=0$ if $i>l(\mu)$. The dual (or conjugate) partition
is denoted by $\widetilde\mu$, i.e.~$\widetilde\mu_j = \max\{ k\,:\,\mu_k\ge j\}$; note 
that $l(\widetilde\mu)=\mu_1$. $\lambda\pm\mu$ denotes the sequence given by 
pointwise addition (resp.~subtraction); while $\lambda+\mu$ is again a 
partition, $\lambda-\mu$ is often not. $(\lambda,\mu)$ denotes the 
concatenation, i.e.~$(\lambda_1,\dots,\lambda_{l(\lambda)},\mu_1,\dots,\mu_{l(\mu)})$. 
Finally $(n^k)$ is the `block' partition $(n,n,\dots,n)$ ($k$ times); and
for $\lambda\subset(n^k)$ we denote its `complement' by $\c\lambda$, 
i.e.~$(\c\lambda)_j=n-\lambda_{k+1-j}$ (we omit the block itself from the 
notation, as it will be always clear from the context).

Let $c_1,c_2,\dots$ and $s_1,s_2,\dots$ be two sequences of (formal) 
variables related by the equation 
\[ (1+c_1t+c_2t^2+c_3t^3+\cdots)\cdot(1-s_1t+s_2t^2-s_3t^3+-\cdots)=1. \]
The \emph{Schur polynomial} $s_\lambda$ is then the determinant 
$s_\lambda = \det[s_{\lambda_i-i+j}] = \det [c_{\widetilde\lambda_i-i+j}]$ (we adopt the convention that
$c_0=1$ and $c_k=0$ for $k<0$; similarly for $s_0$ and $s_{<0}$). Setting
the degree of $c_k$ (resp.~$s_k$) to $k$, $s_\lambda$ will be a homogeneous
polynomial of degree $|\lambda|$. 

If $E\to X$ is a complex vector bundle---or more generally, $E\in K(X)$---, 
we can interpret these sequences as its Chern and Segre classes; the resulting expression is 
denoted by $s_\lambda(E)$, and is an element of $H^*(X)$, which denotes the 
singular cohomology with integer coefficients. 
Note that $s_{\lambda}(-E) = s_{\widetilde\lambda}(E^*)$ 
where $E^*=\Hom(E,\C)$ is the dual bundle of $E$.

The Littlewood-Richardson coefficients will be denoted by $c^\lambda_{\mu,\nu}$, 
i.e.~we have the expansion 
$s_\mu\cdot s_\nu = \sum_\lambda c^\lambda_{\mu,\nu} s_\lambda$.

%--------------------------------------------------------------------
\section{Thom polynomials} 
 
We use the general framework of `Thom polynomials for group actions' introduced by M. Kazarian \cite{kazarianth}; see also \cite{fr}. 

Let $\rho:G\to \GL(V)$ a representation of the Lie group $G$ on the vector space $V$. Then any closed invariant subvariety $\Sigma$ of $V$ represents an equivariant cohomology class $[\Sigma]\in H^*_G(V)\iso H^*(BG)$. We sometimes call this class Thom polynomial 
because $H^*(BG)$ is---at least rationally---a polynomial ring, and Thom polynomials of singularities are special cases where $V$ is the---infinite dimensional---vector space of holomorphic germs $(\C^n,0)\to (\C^p,0)$, 
and $G$ is the group of germs of biholomorphisms $\mathcal{A}(n,p)=j^\infty\Diff(n)\times j^\infty\Diff(p)$. 
Some caution is required in this case since there is no natural topology defined 
on the group $\mathcal{A}(n,p)$, 
so the classifying space $B\mathcal{A}(n,p)$ is not defined. As it is explained 
in \cite{rl} one can work with the classifying space of the subgroup 
of linear germs $\GL_n\times\GL_p\subset\mathcal{A}(n,p)$ instead. 
$H^*(B(\GL_n\times\GL_p))\iso \Z[a_1,\dots,a_n,b_1,\dots,b_p]$ and one 
can interpret the variables as Chern classes: $a_i=c_i(A)$ and $b_i=c_i(B)$ for the 
tautological vector bundles $A^n\to B\GL_n$ and $B^p\to B\GL_p$.

We have a geometric meaning of the Thom polynomial. 
Suppose that $E\to X$ is a $\rho$-bundle, i.e. $E$ is of the form 
$P\times_\rho V$ for some principal $G$-bundle $P\to X$. Then we can 
define a subset $\Sigma(E)$ of the total space of $E$, the union of 
$\Sigma$-points in each fiber.
\begin{prop} 
If $\sigma:X\to E$ is a generic section (transversal to $\Sigma(E)$) 
then for the cohomology class $[\Sigma(\sigma)]\in H^*(X)$ defined by 
\[\Sigma(\sigma):=\{x\in X:\sigma(x)\in \Sigma(E)\}\]
we have 
\[ [\Sigma(\sigma)]=k^*[\Sigma] \]
where $k:X\to BG$ is the classifying map of the principal $G$-bundle $P\to X$.
\end{prop}

%--------------------------------------------------------------------
\section{Thom polynomials of $\Sigma^{i,j}$} 

This approach to $\Sigma^{i,j}$ singularities is fairly standard, 
see eg.~\cite{ronga,avgl,kaza:localized}.

Let $V^n$, $W^p$ be complex vector spaces, and 
$\phi:(V,0)\to (W,0)$ a holomorphic germ; then its derivative 
$\d\phi:T_0 V\to T_0 W$ can be identified with an element of 
$\Hom(V,W)$. Similarly if $f:N\to P$ is a holomorphic map of 
manifolds, then $\d f$ gives a section of $\Hom(TN,f^*TP)$. To define 
the second derivative is more subtle. From the second partial 
derivatives of $\phi$ we can build a linear map in $\Hom(\Sym^2 V, W)$, 
however it is easy to see that for a map $f:N\to P$ we don't get a canonical
section of $\Hom(\Sym^2(TN),f^*TP)$. Using a connection on $TN$ we can get a 
section, but this will depend on the choice of connection. On the other 
hand $f$ defines a section $j^2f$ of the jet bundle $J^2(N,P)$. The space 
$J^2(N,P)$ is a fiber bundle over the product space $N\times P$ (see 
e.g.~\cite{avgl} or \cite{boardman}). Its fibers are diffeomorphic to 
the vector space $J^2(n,p)=\Hom(\C^n,\C^p)\oplus \Hom(\Sym^2\C^n,\C^p)$, 
but the linear structure is not canonical on them. In other words, the 
structure group of $J^2(N,P)$, the Lie group of $2$-jets of germs of 
biholomorphisms
$j^2\Diff(n)\times j^2\Diff(p)$ can be reduced to the
subgroup $\GL_n\times\GL_p$. This reduction is not canonical but 
is a homotopy equivalence, so the equivariant 
cohomology for the two groups are the same.

This phenomenon motivates the definition of the intrinsic quadratic derivative of $f$ (due to I.~R.~Porteous). Let $x\in N$:

\begin{prop} There is a unique linear map 
\[ \d^2_xf:\Sym^2\big(\ker (\d_xf)\big) \to \coker (\d_xf) \]
such that for any $X,Y\in \ker (\d_xf)$ and vector fields $\xi,\eta$ with $\xi_x=X$, 
$\eta_x=Y$ and $g$ holomorphic function defined in a neighborhood of $f(x)$ in $P$:
\[\big(\d^2_xf(X,Y)\big)g\equiv\xi\eta(g\circ f)\;\;\textrm{mod}\;\;\,\im(\d_xf).\]
\end{prop}

The proof is a simple local computation.

We can use the intrinsic second derivative to redefine $\Sigma^{i,j}$. 
First we define the {\em 2-corank} $\corank_2\alpha$ of a map 
$\alpha\in \Hom(\Sym^2V,W)$ where $V,W$ are vector spaces. 
Since $\Sym^2V\cong (V\otimes V)/\Lambda^2V$ there is a canonical inclusion 
$\Hom(\Sym^2V,W)\subset \Hom(V\otimes V,W)$, so $\alpha$ defines an element 
$\tilde\alpha\in \Hom(V,V^*\otimes W)$.

\begin{defn} $\corank_2(\alpha):=\corank(\tilde\alpha)$. \end{defn}

\begin{prop}\label{def} $\Sigma^{i,j}(f)=\{x\in N:\corank (\d_xf)=i,\;\, 
\corank_2 (\d^2_xf)=j$\}. \end{prop}
\noindent
We would like to use this form to show that the Thom polynomial of $\Sigma^{i,j}$ agrees with a Thom polynomial corresponding to a finite dimensional representation.
\begin{defn}\label{defsigma}
Let $V$, $W$ be vector spaces and define the subvarieties
\begin{eqnarray*} 
\Sigma^{\bullet,j}(V,W)&:=&\big\{\alpha\in \Hom(\Sym^2V,W):\corank_2(\alpha)=j\big\} \\
\Sigma^{i,j}(V,W)&:=&\big\{(\alpha_1,\alpha_2)\!\in\!\Hom(V,W)\oplus \Hom(\Sym^2V,W)\,:\corank(\alpha_1)=i, \\
& &\quad \hat\alpha_2\in \Sigma^{\bullet,j}\big(\ker(\alpha_1),\coker(\alpha_1)\big)\big\}
\end{eqnarray*}
where $\hat\alpha_2:\ker(\alpha_1)\to\coker(\alpha_1)$ is the obvious map induced by $\alpha_2$.
\end{defn}
Notice that $\bar\Sigma^{i,j}(V,W)$ is a $G=\GL(V)\times \GL(W)$-invariant subvariety of $\Hom(V,W)\oplus \Hom(\Sym^2V,W)$ so we can define the Thom polynomial $[\Sigma^{i,j}(V,W)]\in H^*(BG)$ (we use the convention that the Thom polynomial of an invariant subset is the equivariant cohomology class defined by its closure).

We will continue with the formalism that $A^n$ resp.~$B^p$ will denote
complex vector spaces equipped with the standard representation of
$\GL(A)\isom\GL_n$ (resp.~$\GL_p$). We will think of them 
either as vector spaces, representations or equivariant vector bundles
over the one-point space; from this last viewpoint they have
(equivariant) Chern classes $a_1,a_2,\dots,a_n$ (resp.~$b_1,\dots,b_p$),
which we will often treat as formal variables. 

As a corollary of the Thom transversality theorem and the definitions above we get:
\begin{prop} 
\[ [\Sigma^{i,j}(n,p)]=[\Sigma^{i,j}(A^n,B^p)]\in H^*\big(B(\GL(A)\times\GL(B))\big)\iso \Z[a_1,\dots,a_n,b_1,\dots,b_p]. \]
\end{prop}

\noindent
As we mentioned before the polynomials $[\Sigma^{i,j}(n,p)]$  
are stable in the following sense:
\begin{prop}\label{stable} 
There exist polynomials $[\Sigma^{i,j}(r)]\in \Z[c_1,c_2,\dots]$ such that for all pairs of natural numbers $(n,p)$ with $p-n=r$ we have 
\[ [\Sigma^{i,j}(n,p)]=[\Sigma^{i,j}(r)](B^p-A^n),\]
where the right hand side means that we substitute $c_k(B^p-A^n)$ into $c_k$ in the polynomial $[\Sigma^{i,j}(r)]$. 
\end{prop}
A more general theorem is due to J.~Damon. 
A proof can be found e.g.~in \cite{fr}. 
Here $B-A$ denotes the difference in the Grothendieck group; $c_i(B-A)$ can 
be interpreted as the $i$'th Taylor coefficient of the formal power series
$(\sum_{k\ge 0} b_kt^k)/(\sum_{l\ge 0} a_lt^l)$. For example
\begin{eqnarray*}
c_1(B-A)&=&b_1-a_1\\
c_2(B-A)&=&b_2-a_2+a_1^2-a_1 b_1
\end{eqnarray*}
and so on. 

%--------------------------------------------------------------------
\section{Calculation of the Thom polynomials}

It follows from 
Definition \ref{defsigma}
that $\Sigma^{i,j}(A^n,B^p)$ is empty if $n<i$ and for $n=i$ it has a particularly simple structure:
\[\Sigma^{i,j}(A^i,B^{r+i})=\big\{(0,\alpha_2)\,:\,\alpha_2\in\Sigma^{\bullet,j}(A^i,B^{r+i})\big\},\]
so its Thom polynomial is a product:
\[ \big[\Sigma^{i,j}(A^i,B^{r+i})\big]=e\big(\Hom(A^i,B^{r+i})\big)\!\cdot\!\big[\Sigma^{\bullet,j}(A^i,B^{r+i})\big],\]
where $e(\Hom(A^i,B^{r+i}))$ denotes the Thom polynomial of 
$\{0\}\subset \Hom(\C^i,\C^{r+i})$---the (equivariant) \emph{Euler class} 
of this representation---and $[\Sigma^{\bullet,j}(A^i,B^{r+i})]$ is the Thom 
polynomial of the subvariety 
$\bar\Sigma^{\bullet,j}(A^i,B^{r+i})\subset \Hom(\Sym^2 A^i,B^{r+i})$. 
Since the Euler class of a representation is the product of its weights 
we get that $e=\prod(\beta_j-\alpha_i)$ where $\alpha_i$ and $\beta_j$ are 
the {\em Chern roots} of 
$A$ and $B$, i.e.~$a_k$ is the $k\textsuperscript{th}$ elementary 
symmetric polynomial of the variables $\alpha_i$ and similarly $b_k$ 
is the $k\textsuperscript{th}$ elementary symmetric polynomial of the 
variables $\beta_j$. 
Comparing this with Proposition \ref{stable} we get the equations:
{\setlength\arraycolsep{0pt}
\begin{eqnarray}
\label{restreq1}&\big[\Sigma^{i,j}(r)\big]&(B^{r+i-1}-A^{i-1}) \,=\, 0; \\
\label{restreq2}&\big[\Sigma^{i,j}(r)\big]&(B^{r+i}-A^i) \,=\, e\big(\Hom(A,B)\big)\!\cdot\!\big[\Sigma^{\bullet,j}(A,B)\big].
\end{eqnarray}} 

\noindent
Equations (\ref{restreq1}) and (\ref{restreq2}) can be also 
interpreted as \emph{restriction equations} in the sense of \cite{fr}, 
the right hand side of equation (\ref{restreq2}) being an `incidence class' 
in the sense of \cite{rrthompol}. From this point of view equation 
(\ref{restreq1}) follows from the obvious fact that $\Sigma^{i,j}(f)\subset \Sigma^{i}(f)$ for 
any holomorphic map $f$, and equation (\ref{restreq2}) is a consequence
of the local behavior of the set $\Sigma^{i,j}(f)$ at a point 
of $\Sigma^{i}(f)$. This was indeed our first approach and we only later realized 
that (\ref{restreq1}) and (\ref{restreq2}) also follows from stability.\\

The idea of our calculation is to solve the system of equations
(\ref{restreq1}) and (\ref{restreq2}).
To do that, we have 
to study the homomorphism 
\[\rho_{n,p}: \Z[c_1,c_2,\dots] \to \Z[a_1,\dots,a_n,b_1,\dots,b_p]\]
sending $c_i$ to $c_i(B^p-A^n)$. 
Elements in the image of $\rho_{n,p}$ are called 
{\em supersymmetric polynomials} 
(or \emph{Schur functions in difference of alphabets}). The following 
proposition states some of the fundamental properties of supersymmetric 
polynomials (in other words of the map $\rho_{n,p}$).

\begin{prop}\label{factorization}\mbox{}
\begin{enumerate}
\item[(i)\,] $\ker (\rho_{n-1,p-1})=\langle s_\lambda:n^p\subset \lambda\rangle$, where $\langle\ \rangle$ means the generated $\Z$-module.
\item[(ii)] Suppose that $\lambda$ is a partition containing the `block' partition $n^p$. If we also assume that $(n+1)^{(p+1)}\not\subset\lambda$ 
then $\lambda$ is of the form $(n^p+\beta,\alpha)$, where $l(\beta)\le p$ and $\alpha_1\le n$; 
i.e.~$\lambda_i=n+\beta_i$ for $i\leq p$ and $\lambda_i=\alpha_{i-p}$ for $i>p$.
\[ s_\lambda(B^p-A^n) = \left\{ \begin{array}{lll}
(-1)^{|\alpha|}e(\Hom(A,B))s_{\widetilde\alpha}(A)s_\beta(B) &\;& (n+1)^{(p+1)}\not\subset\lambda \\
0 &\;& (n+1)^{(p+1)}\subset\lambda 
\end{array}\right.{}.\]
\end{enumerate}
\end{prop}
\noindent
The proof can be found in any introduction on supersymmetric polynomials,
e.g. in \cite[\S 3.2]{fulton-pragacz}. Part (i) is a corollary of a result 
of Pragacz \cite{pragacz88} on universally supported classes 
({\em avoiding ideal} in the terminology of \cite{fr}) for $\Sigma^i$; 
part (ii) is sometimes called the {\em factorization formula}.

From this formula it is clear that the system of equations above do not
have a unique solution: If we write the solution as a linear combination
of Schur polynomials, we will have an ambiguity in the terms $s_\lambda$
where $(i+1)^{r+i+1}\subset\lambda$. But it is also clear that all the
other terms are uniquely determined by
our equations. To our surprise, these ambiguous terms turned out to be zero:

%\newtheorem*{thm_rongacor}{Theorem \ref{rongacor}} 
%\begin{thm_rongacor}
\begin{thm}\label{rongacor}
Write the universal polynomial $[\Sigma^{i,j}(r)]$ as a linear
combination of Schur polynomials: $[\Sigma^{i,j}(r)]=\sum e^\lambda s_\lambda$,
where $e^\lambda$ are (integer) coefficients. 
Then $e^\lambda=0$ if $(i+1)^{(r+i+1)} \subset \lambda$.   
\end{thm}
%\end{thm_rongacor}

\noindent
The proof, which is based on Ronga's \cite{ronga} pushforward formula for $[\Sigma^{i,j}(r)]$,
is given in Section \ref{review}, Theorem \ref{rongacor_gen}. The same proof yields some more vanishing 
results, but we don't need them.
\\

Theorem \ref{rongacor}, Proposition \ref{factorization} 
and the two equations above together
imply the following:
\begin{corr}\label{connection} 
Write the polynomial 
$\big[\Sigma^{\bullet,j}(A^i,B^{r+i})\big]$ 
as a linear combination of product of Schur polynomials in the 
variables $a_i$ and $b_i$:
\[\big[\Sigma^{\bullet,j}(A^i,B^{r+i})\big]=\sum e_{\alpha,\beta} s_\alpha(A)s_\beta(B).\]
Then 
\[ [\Sigma^{i,j}(r)]=\sum (-1)^{|\alpha|}e_{\alpha,\beta}s_{(i^{(r+i)}+\beta,\widetilde\alpha)}.\]
\end{corr}

 Calculating $[\Sigma^{\bullet,j}(A^i,B^{r+i})]$ \emph{in this form}
seems to be a very difficult problem in general; although the difficulties
are purely combinatorial, as we have the following pushforward formula (see
also Lemma \ref{gysin} about calculation of pushforward).

\begin{prop}
Let $\pi:\Gr_j(A^i)\to\pt$ denote the projection map from the Grassmannian
of $j$-planes in $A$ to the one-point space, and
$0\to R^j\to \pi^*\!A\to Q^{i-j}\to 0$ the tautological exact sequence (of equivariant
vector bundles) over $\Gr_j(A)$. Then 
\[\big[\Sigma^{\bullet,j}(A,B)\big]=
\pi_* c_{\mathrm{top}}\!\left(\pi^*\!B\otimes(R\otimes Q + \Sym^2R)^*\right). \]
\end{prop}

Let us emphasize that in the light of Corollary \ref{connection} above, this formula,
while much simpler than Ronga's or Kazarian's pushforward formulas
for $[\Sigma^{i,j}(r)]$, contains the same amount of information. The proof of the
proposition is analogous to what is written in the first part of Section 3 in \cite{lpschurq};
we don't repeat it here, as we won't use this formula in the rest of the paper.\\

There are two particular cases when we know the solution in the required form, 
namely, the cases $r+i=1$ and $j=1$. It is not hard too see why these are simpler from the pushforward
formula above.

\begin{thm}\label{sigma_ht} 
$[\Sigma^{\bullet,j}(A^i,L^1)]=2^j\!\cdot\!s_\delta(A^*\otimes \sqrt{L})$ 
where $\delta=(j,j-1,\dots,2,1)$.
\end{thm}

\noindent
Note that the line bundle $L$ has no square root, so the formula 
above should be understood formally: the only Chern root of 
$\sqrt{L}$
is $\beta/2$ where $\beta=\beta_1$ is the Chern root of $L$, 
and then the Chern roots of $A^*\otimes \sqrt{L}$ are 
$-\alpha_1+\beta/2,\dots,-\alpha_n+\beta/2$.

\begin{proof} Notice that the elements of $\Hom(\Sym^2\C^i,\C)$ can be 
identified with symmetric $i\times i$ matrices and then $\corank_2=\corank$,
so the Thom polynomial in question is given by the {\em twisted symmetric 
degeneracy locus} formula (\cite{harris-tu}, \cite{jlp}, \cite{pragacz2}, \cite{fulton2}).
This formula first appeared in \cite{harris-tu}, but their proof was 
erroneous; the correct proof was given in \cite{fulton2}.
A general explanation of twisting can be found in \cite{fnr_forms}.
\end{proof}

\begin{thm} \label{sigma_porteous}$[\Sigma^{\bullet,1}(A^i,B^{r+i})] = c_{i(r+i-1)+1}(A^*\otimes B - A)$.
\end{thm}

\begin{proof}
The codimension
of $\Sigma^{\bullet,1}(V^n,W^p)\subset \Hom(\Sym^2V,W)$ is $pn-n+1$, which equals to the
codimension of $\Sigma^1(V,V^*\otimes W)\subset \Hom(V,V^*\otimes W)=\Hom(V\otimes V,W)$; so---exactly as noted in \cite{lpschurq}, where a similar degeneracy locus problem is considered---we are in the situation of the Giambelli-Thom-Porteous formula: 
\[ \big[\Sigma^{\bullet,1}(V^n,W^p)\big]=[\Sigma^1(V,V^*\otimes W)]=c_{pn-n+1}(V^*\otimes W-V).\]
\end{proof}

According to Corollary \ref{connection}, the only thing we need 
is the separation of variables in the formulas of Theorem \ref{sigma_ht} 
and \ref{sigma_porteous}. We will use the following lemma, due to Lascoux.

\begin{lem}[\cite{lascoux}]\label{tensor-rule} 
Denote by $E_{\lambda/\mu}(n)$ the determinant
\[ E_{\lambda/\mu}(n) = 
\det \left[ \binom{\lambda_i+n-i}{\mu_j+n-j} \right ]_{i,j\in n\times n}. 
%= \det \left[ \binom{\lambda_i+n-i}{\lambda_i - \mu_j - i + j} \right ]_{i,j\in n\times n}.
\]
\begin{enumerate}
\item Let $A^n$ and $B^p$ be an $n$-dimensional and a $p$-dimensional 
vector bundle, respectively. Then
\[ \sum_{k} c_k(A\otimes B) = 
\sum_{\mu\subset\lambda\subset p^n} E_{\lambda/\mu}(n) s_\mu(A) s_{\c\widetilde\lambda}(B) 
%= \sum_{\mu\subset\lambda\subset n^p} E_{\lambda/\mu}(p) s_{\c\widetilde\lambda}(A)  s_\mu(B) 
\]
\item Furthermore, if $L$ is a line bundle and $\lambda$ is partition with $l(\lambda)\le n$, then
\[ s_\lambda(A\otimes L) = 
\sum_{\mu\subset\lambda}E_{\lambda/\mu}(n)\!\cdot\! c_1(L)^{|\lambda|-|\mu|}\!\cdot\!s_\mu(A) 
\]
\end{enumerate}
\end{lem}

\begin{rem}
The coefficients $E_{\lambda/\mu}(n)$ are known to be nonnegative.
This is for example a very special case of \cite{pragaczpar}, Corollary 7.2 which says
that if we set 
\[ s_\lambda(S^{\mu_1}E_1\otimes S^{\mu_2}E_2 \otimes\cdots\otimes S^{\mu_k}E_k) = 
\sum E_{\lambda,\blds{\mu}}^{\blds{\nu}} s_{\nu_1}(E_1)s_{\nu_2}(E_2)\cdots s_{\nu_k}(E_k), \]
where $S^\mu$ is the Schur functor associated to the partition $\mu$, 
then all the coefficients $E_{\lambda,\blds{\mu}}^{\blds{\nu}}$ will be nonnegative.
\end{rem}

\noindent
The lemma, together with Theorem \ref{sigma_ht} and the framework built above immediately implies one of our main theorems:

\begin{thm}\label{main1}
The Thom polynomial of the second order Thom-Boardman singularity $\Sigma^{i,j}(-i+1)$ is 
\[ \big[\Sigma^{i,j}(-i+1)\big] =  \sum_{\mu\subset\delta} 2^{|\mu|-j(j-1)/2} \cdot E_{\delta/\mu}(i) \cdot s_{(d-|\mu|,\widetilde\mu)} \] 
where $\delta$ is the `staircase' partition $\delta=(j,j-1,\dots,2,1)$ and  
\[ d=\codim\,\Sigma^{i,j}(-i+1)=i+|\delta|=i+{\textstyle \binom{j+1}{2}}.\]
\end{thm}

\noindent
Similarly, Theorem \ref{sigma_porteous} leads to 

\begin{thm}\label{main2} Using the shorthand notation $h=r+i$,
\[{}\big[\Sigma^{i,1}(r)\big] = 
\sum_{(\lambda,\mu)\in J} s_{(i^{h}+\lambda,\mu)}
\cdot \!\!\! \sum_{x\in\{0,1\}^{l(\mu)}} \!\!\! E_{\c \widetilde\lambda / (\mu-x)\phtilde}(i) \]
where $J=\{(\lambda,\mu)\::\:\lambda\subset i^{h},\; \mu_1\le i,\; |\lambda|+|\mu|=ih-i+1,
\textrm{ and $\mu-x$ is a valid partition}\}$.
\end{thm}

\begin{proof}
According to Corollary \ref{connection}, to solve the equation for
$[\Sigma^{i,1}]$ all we have to do is to expand the formula
$c_{ih-i+1}(A^*\otimes B - A)$ into linear combination of products
of Schur polynomials. For the sake of convenience, we calculate the
total Chern class
\[ \sum_{m\ge 0} c_m(A^*\otimes B - A) = 
\left(\sum_{k\ge 0} c_k(A^*\otimes B) \right) 
\cdot\left(\sum_{l\ge 0}c_l(-A)\right). \]
Using Lemma \ref{tensor-rule}, the Pieri formula, and 
\[ c(-A) = \sum_{l\ge 0} c_l(-A) = \sum_{k\ge 0} (-1)^i s_k(A), \]
we will get
\[ c(A^*\otimes B - A) = \sum_{\mu\subset\lambda\subset i^h}
\sum_{x\in\{0,1\}^{l(\mu)}} (-1)^{|\mu+x|} E_{\widetilde\lambda/\widetilde\mu}(i)\cdot
s_{(\mu+x)\phtilde}(A) s_{\c\lambda}(B), 
\]
where the second sum runs over $0$-$1$ sequences such that $\mu+x$ is
a valid partition.
From this the theorem follows directly, using 
the fact that 
$E_{\lambda/\mu}(k)=0$ if $\mu\not\subset\lambda$ and $k\ge l(\lambda),l(\mu)$.
\end{proof}

Note that in both cases, the Thom polynomial is a nonnegative linear 
combination of Schur polynomials. Based on computational evidence,
we can formulate the following

\begin{conj}The Thom polynomials of all Thom-Boardman classes can be 
written as nonnegative linear combination of Schur polynomials.
\end{conj}

With some work, we can get an even more spectacular formula. 
Introduce the shorthand notations
\[ \gbinom{n}{k} := \sum_{j=0}^k \binom{n}{j}
\quad\quad\textrm{and}\quad\quad
 F_{\lambda/\mu}(n):=
\det \left[ \gbinom{\lambda_k+n-k}{\mu_l+n-l} \right ]_{k,l\in n\times n}.
\]
Note that the numbers $\gbinom{n}{k}$ also form a Pascal-like triangle: 
\[\setlength{\arraycolsep}{2pt}
\begin{array}{*{11}{p{8pt}}}
&&&&& 1 &&&&& \\
&&&& 1 && 2 &&&& \\
&&& 1 && 3 && 4 &&& \\
&& 1 && 4 && 7 && 8 && \\
& 1 && 5 && 11 && 15 && 16 & \\
1 && 6 && 16 && 26 && 31 && 32 
\end{array}\]

\begin{thm}\label{main2nice}
\[ {}[\Sigma^{i,1}(r)] = 
\sum_{(\nu,\mu)\in I}
F_{\nu/\mu}(i) \cdot s_{(i^h+\c\widetilde\nu,\widetilde\mu)} 
\]
where $I=\{(\nu,\mu)\::\: \nu\subset h^i,\; l(\mu)\le i,\textrm{ and }\;|\nu|-|\mu|=i-1\}$.
\end{thm}

\begin{rem}
Note that the coefficients do not depend on the relative codimension $r=h-i$. 
This is not a coincidence, 
and a similar phenomenon occurs for a large class of singularities; see 
\cite{fr-d-stab} and the discussion in the introduction.
\end{rem}

\begin{proof}
According to Theorem \ref{main2}, the coefficient $a_{\nu,\mu}$ of $s_{(i^h+\c\widetilde\nu,\widetilde\mu)}$ is
a sum, which we can rewrite as follows:
\[ a_{\nu,\mu} = \sum_{x\in\{0,1\}^{\mu_1}} E_{\nu / (\widetilde\mu -x)\phtilde}(i) = 
\sum_{\alpha_1 = \mu_2}^{\mu_1} \sum_{\alpha_2 = \mu_3}^{\mu_2} \cdots
\sum_{\alpha_i = 0}^{\mu_i} E_{\nu / \alpha}(i) \]
Expanding the determinant $E_{\nu/\alpha}(i)$ and rearranging the sums yields
\[ a_{\nu,\mu} = \det\left[
\gbinom{\nu_k + i - k}{\mu_l + i - l} - \gbinom{\nu_k+i-k}{\mu_{(l+1)} + i - (l+1)}
\right]_{k,l\in i\times i} \]
Observe that $a_{\nu,\mu}$ is of the form $\det(A-B)$ where
\[ B_{k,l} = \left\{\begin{array}{lll} 
A_{k,l+1} &\quad & \textrm{if }l<n \\ 
0         &\quad & \textrm{if }l=n
\end{array}\right.
\]
It is an easy exercise then to prove that in such a situation 
$ \det(A+\beta B)=\det(A)$
holds for any $\beta\in\C$.
\end{proof}

%--------------------------------------------------------------------
\section{Review of Ronga's formula} \label{review}

Ronga's result \cite{ronga}, expressed in our language, is the following. 
Let $X^n\to M$ be a complex vector bundle; 
$p:\Gr_i(X)\to M$ the bundle of Grassmannians of $i$-planes in $X$, 
$Y^i\to \Gr_i(X)$ the tautological subbundle over $\Gr_i(X)$; finally,
$\pi:\Gr_j(Y)\to \Gr_i(X)$ the bundle of Grassmannians of $j$-planes in $Y$
and $0\to R^j\to\pi^*Y\to Q^{i-j}\to 0$ the tautological exact sequence
of vector bundles over $\Gr_j(Y)$. Furthermore, introduce the 
shorthand notations $h=r+i$ and $k=ij-\binom{j}{2}$. Now
\[ \big[\Sigma^{i,j}(r)\big](-X) = 
(-1)^{hk} p_*\pi_*\Big[
s_{(n+r)^i}(\pi^*Y^*)\cdot s_{(h^k)}\big( R\otimes Q + \Sym^2R + \pi^*p^*X - \pi^*Y\big)
\Big] \]
where the left hand side means that we substitute the Chern classes of $-X$
(which are the same as the Segre classes of $X^*$)
into the universal polynomial $[\Sigma^{i,j}(r)]$.
This formula is not very well-suited for direct computations; nevertheless,
we can use it to get some qualitative information about these Thom polynomials.

\begin{thm}\label{rongacor_gen}     
Write the universal polynomial $[\Sigma^{i,j}(r)]$ as a linear
combination of Schur polynomials: $[\Sigma^{i,j}(r)]=\sum e^\lambda s_\lambda$,
where $e^\lambda$ are (integer) coefficients. Then $e^\lambda=0$ if $\lambda$
satisfies any of the following three conditions:
\begin{enumerate}
\item[(a)] $i^{h} \not\subset \lambda$     
\item[(b)] $(i+1)^{(h+1)} \subset \lambda$ 
\item[(c)] $\lambda_1 > i+k$. 
\end{enumerate}
\end{thm}

\noindent
We will use the following well-known lemma (see e.g.~\cite{pragacz88})
about pushforwards (or Gysin maps):

\begin{lem}\label{gysin}
Let $E^n\to M$ be a complex vector bundle, $\pi:\Gr_r(E)\to M$ the bundle
of Grassmannians of $r$-planes in $E$, and $0\to R^r \to \pi^*E \to Q^q \to 0$
the tautological exact sequence of bundles over $\Gr_r(E)$. Then
\[\pi_*\big[s_\mu(R)s_\nu(Q)\big] = s_{(\nu-r^q,\mu)}(E).\]
\end{lem}

\begin{rem}
This formula should be understood as follows: $(\nu-r^q,\mu)$ is very often
not a valid partition; but we can extend the definition of the Schur
polynomials for arbitrary sequences. Every such ``generalized Schur polynomial''
agrees with zero or a ``usual Schur polynomial'' up to sign. For example for 
the particular case $\nu=0$ the formula gives
\[ \pi_*s_\mu(R) = s_{(-r^q,\mu)}(E) = (-1)^{qr}s_{(\mu-q^r)}(E), \]
which is zero if $q^r\not\subset\mu$.
Note also that $\pi_*$ is a $H^*(M)$-module homomorphism.
\end{rem}

\begin{proof}[Proof of Theorem \ref{rongacor_gen}]
With some abuse of notation, we will omit the pullbacks from the formulas;
that is, we will simply write $Y$ instead of $\pi^*Y$ and so on.
All three 
claims will be the consequence of the following computation.
First, using the expansion 
\[ s_{\lambda}(A+B)=\sum_{\mu,\nu} c^\lambda_{\mu,\nu}\cdot s_{\mu}(A)s_{\nu}(B) \]
---which, for the special case $\lambda=h^k$ gives
$s_{h^k}(A+B) = \sum_{\mu\subset h^k}s_\mu(A)s_{\c\mu}(B)$---
we get
\[ \big[\Sigma^{i,j}(r)\big](-X) = 
(-1)^{hk}\sum_{\lambda\subset h^k} 
p_*\Big[ s_{(n+r)^i}(Y^*) s_\lambda(X-Y) \cdot \pi_* s_{\c\lambda}(R\otimes Q + \Sym^2R) 
\Big]. \]
We are not interested in the exact result of the inner pushforward; instead
we just set
\[  (-1)^{hk}\!\cdot\!\pi_* s_{\c\lambda}(R\otimes Q + \Sym^2R) = \sum_{l(\mu)\le i} f^\mu_\lambda\cdot s_\mu(Y), \]
where $f^\mu_\lambda$ are some coefficients.
Using the above expansion again, now for $s_\lambda(X-Y)$ we get:
\[ \big[\Sigma^{i,j}(r)\big](-X) =  
\sum_{\lambda\subset h^k} \sum_{\alpha,\beta\subset\lambda} \sum_{l(\mu)\le i}
c^\lambda_{\alpha,\beta}f^\mu_\lambda \cdot s_\alpha(X) 
\cdot p_*\Big[ s_{(n+r)^i}(Y^*) s_{\widetilde\beta}(Y^*) s_\mu(Y) 
\Big]. \]
Using the Littlewood-Richardson rule, Lemma \ref{gysin} and that 
the rank of $Y$ is $i$,
it follows immediately that
\[ p_*\Big[ s_{(n+r)^i}(Y^*) s_{\widetilde\beta}(Y^*) s_\mu(Y) 
\Big] = \sum_{l(\gamma)\le i} g_\gamma \cdot s_{(h^i+\gamma)}(X), \]
where the $g_\gamma$'s are integer coefficients. Now, we see that 
$[\Sigma^{i,j}](-X)$ is a linear combination of terms of the
form $s_\alpha(X)s_{(h^i+\gamma)}(X)$, where $\alpha\subset h^k$
and $l(\gamma)\le i$. From the Littlewood-Richardson rule it
follows directly that the expansion of such a term satisfies 
the \emph{duals} of all three claims of the theorem, that is, the duals of the partitions
appearing in the expansions satisfies the three conditions; 
thus, using the identity $s_\lambda(-X)=s_{\widetilde\lambda}(X^*)=(-1)^{|\lambda|}s_{\widetilde\lambda}(X)$
the theorem follows.
\end{proof}

%--------------------------------------------------------------------
\section{Examples} 

\begin{exmp} 
The Thom polynomials of the singularities $\Sigma^{i,j}(-i+1)$ for $j\le 2$ are 
\begin{eqnarray*} 
\big[\Sigma^{i,0}\big] & = & s_{i} \\ [1pt]
\big[\Sigma^{i,1}\big] & = & is_{i+1} + 2s_{i,1}  \\[1pt]
\big[\Sigma^{i,2}\big] & = & \textstyle\binom{i+1}{3}s_{i+3} + (i^2-1)s_{i+2,1} + 2(i+1)s_{i+1,2} + 2(i-1)s_{i+1,1,1} + 4s_{i,2,1}
\end{eqnarray*} 
\end{exmp} 

\theoremstyle{definition} 
\newtheorem*{ex_morin}{Morin singularities}  
\begin{ex_morin}
Recall that the Morin singularity $A_2(r)$, where $r$ is the relative
codimension $r=p-n$, is $A_2(r)=\Sigma^{1,1}(r)$ for $r$ nonnegative 
and $A_2(r)=\Sigma^{1-r,1}(r)$  for $r$ negative. We have
\[ \big[A_2(r)\big] = \left\{\begin{array}{lll} \displaystyle
\sum_{k=0}^{r+1} 2^k s_{2^{r+1-k},1^{2k}} &\;\;&\textrm{if } r\ge 0 \\ 
2s_{1-r,1} +(1-r)s_{2-r} &\;\;&\textrm{if } r\le 0 \phantom{\bigg|} \end{array}\right.\]
The $r\ge 0$ case was already known, see \cite{ronga}.
\end{ex_morin} 

\newtheorem*{ex_sigma21}{Thom polynomials of $\Sigma^{2,1}$}  
\begin{ex_sigma21}Let $h$ denote $h=r+2$. Then we have
\[ {}[\Sigma^{2,1}(r)] =  
\sum_K
\left( \gbinom{a+1}{d+1}\gbinom{b}{c} - \gbinom{a+1}{c}\gbinom{b}{d+1}\right )
\cdot s_{(h+d,h+c,h-b,h-a)\phtilde}
\]
where $K=\{(a,b,c,d)\in\N^4\::\:b\le a \le h,\; c\le d,\; c+d=a+b-1\}$.
\end{ex_sigma21}

\begin{rem} 
We can state the analogous theorems for real singularities 
with cohomology in $\Z_2$-coefficients by replacing Chern classes with 
the corresponding Stiefel-Whitney classes. The Thom polynomials for the \emph{real} 
$\Sigma^{i,1}(-i+1)$ were already calculated by Thom in \cite{thom}.
\end{rem}

\bibliography{sigmaij} 
\bibliographystyle{amsalpha} 
 
\end{document}